\documentclass[12pt,a4paper,twoside]{article}

\usepackage[centertags]{amsmath}
\usepackage[T1]{fontenc}
\usepackage{amsfonts}
\usepackage{amssymb}
\usepackage{ntheorem}
\usepackage{ulem}
\usepackage{epsfig}
\usepackage{subfigure}
\usepackage[all]{xy}
\usepackage{changebar}
\usepackage{color}
\usepackage{mathrsfs}
\usepackage[boldsans]{ccfonts}
\usepackage{array}
\usepackage{tabulary}
\usepackage{supertabular}
\usepackage{calc}
\usepackage{eepic}
\usepackage{titlesec}
\usepackage{graphicx}
\usepackage{enumitem}

\theoremstyle{break} \newtheorem{theorem}{Theorem}[section]
\theoremstyle{break} 
\theoremstyle{break} 
\theoremstyle{nonumberbreak}  
\theoremstyle{break} 
\theoremstyle{break} \newtheorem{corollary}[theorem]{Corollary}
\theoremstyle{break} 
\theoremstyle{break} 
{\theorembodyfont{\rmfamily}
\newtheorem{remark}[theorem]{Remark}
{\theorembodyfont{\rmfamily}}
\theoremstyle{break} 
\theoremstyle{break} 
\theoremstyle{break} 
\theoremstyle{break} 
\numberwithin{equation}{section}

\newcommand{\hide}[1]{}

\newcommand{\R}{{\mathbb{R}}}
\newcommand{\D}{{\mathbb{D}}}
\newcommand{\Ha}{{\mathbb{H}}}
\newcommand{\C}{{\mathbb{C}}}

\renewcommand{\H}{{\mathbb{H}}}
\def\dD{\mathop{{\rm d}_\D}}

\def\Re{\mathop{{\rm Re}}}
\def\Im{\mathop{{\rm Im}}}

\def\arcsinh{\mathop{\rm Arcsinh}}
\def\Shypp{\mathop{\text{S}^{+}_{\rm h}}}
\def\Shypm{\mathop{\text{S}^{-}_{\rm h}}}
\def\Shyppm{\mathop{\text{S}^{\pm}_{\rm h}}}
\def\Seucp{\mathop{\text{S}^{+}_{\rm e}}}
\def\Seucm{\mathop{\text{S}^{-}_{\rm e}}}

\begin{document}

\begin{center}
{\Large \bf Rogosinski's lemma for univalent functions, \\[3mm] hyperbolic
  Archimedean spirals and the Loewner equation}

\end{center}
\renewcommand{\thefootnote}{\arabic{footnote}}
\setcounter{footnote}{0}

\smallskip

\begin{center}
{\large Oliver Roth and Sebastian Schlei{\ss}inger}\\
\today
\\[3mm]
\end{center}

\smallskip

\begin{center}
\begin{minipage}{13.2cm}
{\bf Abstract.}{ \small
We describe the region $\mathcal{V}(z_0)$ of values of $f(z_0)$ for all normalized bounded
univalent functions $f$ in the unit disk $\D$ at a fixed point $z_0 \in \D$.
The proof is based on 
 the radial Loewner differential equation. We also prove an analogous result for
 the upper half-plane using the chordal Loewner equation.}
 \end{minipage}
\end{center}

\smallskip

\section{Results: the unit disk} \label{sec:1}

We denote by
 $\D:=\{z \in \C \, : \, |z|<1\}$  the open unit disk in the
complex plane $\C$ and by $\mathcal{H}_0(\D)$  the
 set of all holomorphic functions $f : \D \to \D$ normalized by $f(0)=0$ and
 $f'(0) \ge 0$. The Schwarz lemma tells us that $|f(z_0)| \le |z_0|$ for any
 $f \in \mathcal{H}_0(\D)$ and any $z_0 \in \D$.
In 1934 Rogosinski \cite{Rog1934} (also \cite[p.~200]{Dur1983}) proved a far reaching sharpening of the Schwarz lemma by 
giving an explicit description of the region of values of
$f(z_0)$, $f \in \mathcal{H}_0(\D)$, at a fixed point $z_0 \in \D$.
In this note we prove an analogue of Rogosinski's result 
for \textit{univalent} functions in $\mathcal{H}_0(\D)$ by 
providing an explicit description of the regions of values
$$ \mathcal{V}(z_0):=\left\{ f(z_0) \, : \, f \in \mathcal{H}_0(\D) \text{ univalent}
\right\} \, , \qquad z_0 \in \D\, . $$


It turns out that the set $\mathcal{V}(z_0)$ admits a fairly
appealing description in terms of hyperbolic geometry.
In order to state the main results, we 
 therefore endow the unit disk  with the standard hyperbolic metric
$$\lambda_{\D}(z)\, dz=\frac{2\,|dz|}{1-|z|^2} $$ of constant curvature
$-1$, see \cite{BM2007}.
 The  induced hyperbolic distance $\dD(z,w)$ between $z,w \in \D$ is then given by
\begin{equation} \label{eq:hypdist}
\dD(z,w)=\log \frac{1+\displaystyle 
\left|\frac{z-w}{1-\overline{w} z}
   \right|}{1-\displaystyle 
\left| \frac{z-w}{1-\overline{w} z} \right|} \, .
\end{equation}
\begin{theorem} \label{thm:1}
Let $z_0 \in \D$. Then 
$$ \mathcal{V}(z_0) \cup \{0\}=\left\{z=|z| e^{i \varphi} \in \D \, : \,
\dD(0,z)-\dD(0,z_0)\le -|\varphi-\arg z_0| , \, \varphi \in \R \right\} \, .$$

\begin{figure}[h]
\centerline{\includegraphics[width=4cm]{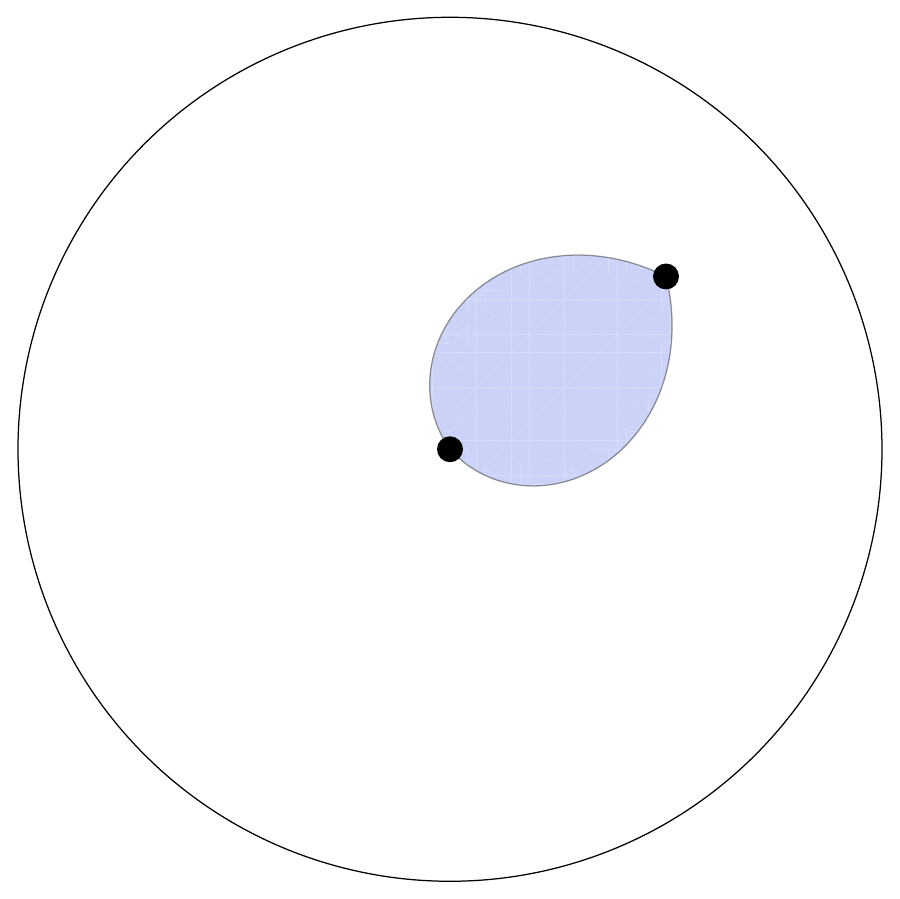} \hspace{2cm}
\includegraphics[width=4cm]{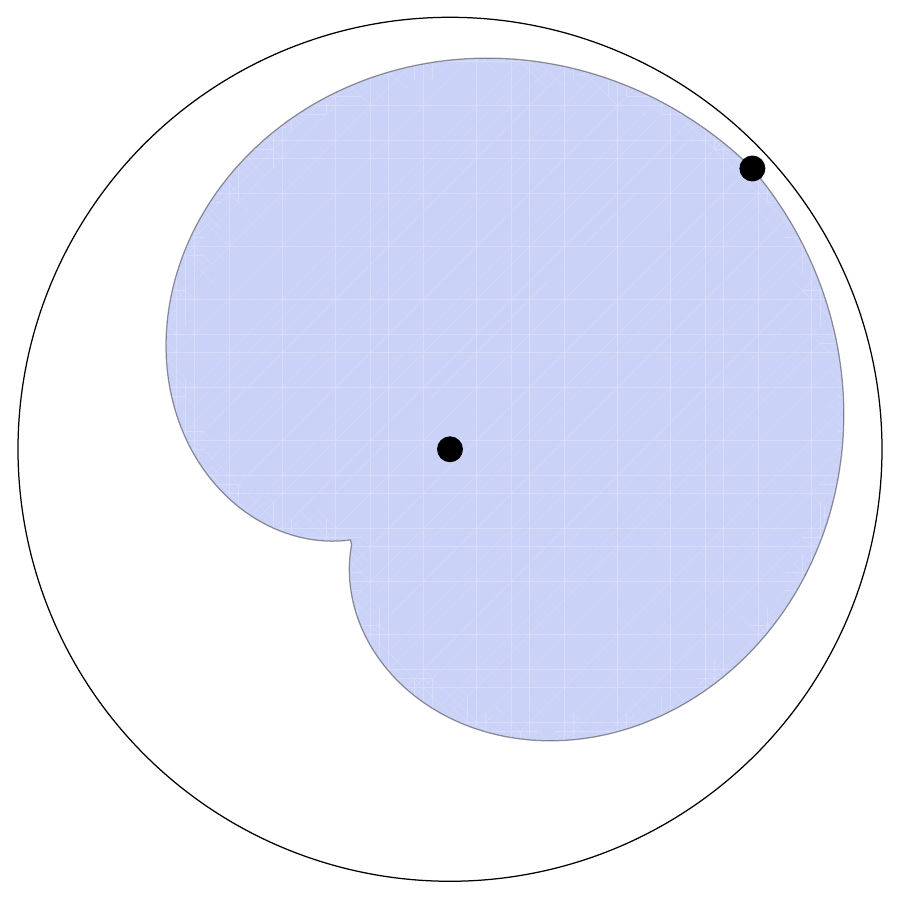}}
\caption{The set $\textcolor{blue}{\mathcal{V}(z_0)}$ for $z_0=0.5+0.4 i$ (left) $z_0=0.7+0.65 i$ (right)
}
\end{figure}
\end{theorem}

In particular, 
the boundary of the region of values $\mathcal{V}(z_0)$ is composed of parts
of the two curves
\begin{eqnarray*}
& &\Shypp(z_0) \, :=\left\{ |z| e^{i \varphi} \, : \,
\dD(0,z)=+\varphi+\dD(0,z_0)-\arg z_0 \, , \, \varphi 
\ge \arg z_0-\dD(0,z_0)\right\}\, ,\\[1mm]
& & \Shypm(z_0) \, :=\left\{ |z| e^{i \varphi} \, : \,
\dD(0,z)=-\varphi+\dD(0,z_0)+\arg z_0 \, , \, \varphi \le \arg
z_0+\dD(0,z_0)\right\} \, .
\end{eqnarray*}

\begin{figure}[h]
\centerline{\includegraphics[width=4cm]{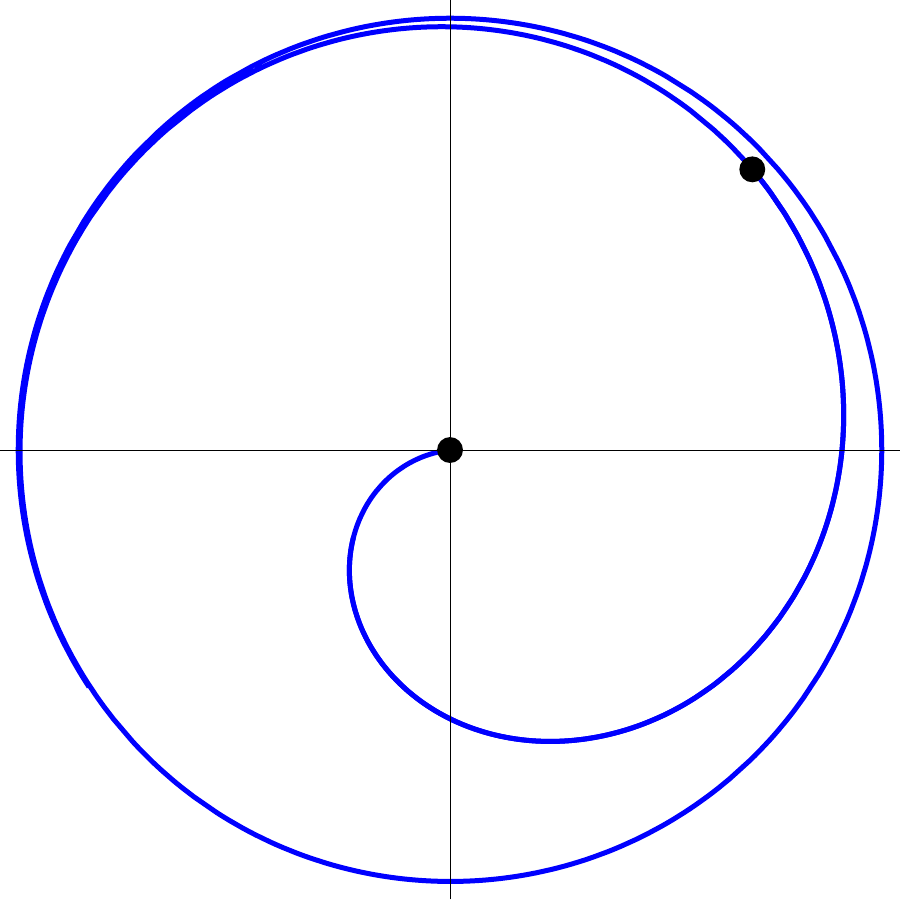}\hspace{2cm}
\includegraphics[width=4cm]{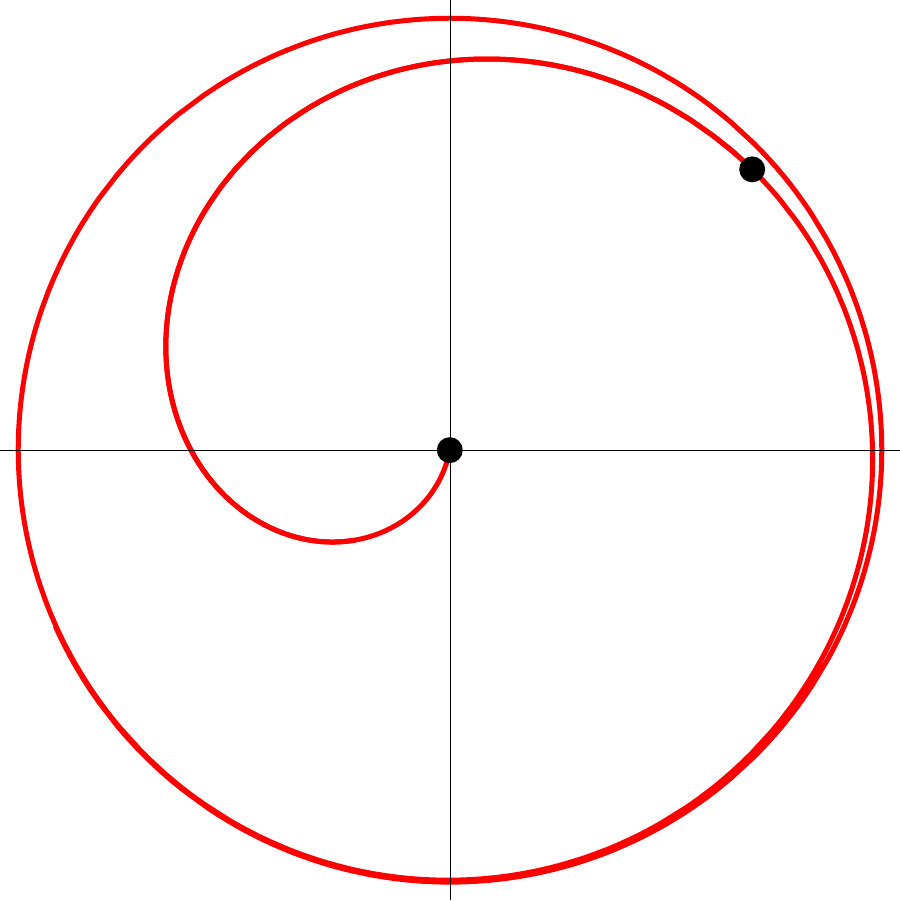}}
\caption{Hyperbolic Archimedean spirals $\textcolor{blue}{\Shypp(z_0)}$ and $\textcolor{red}{\Shypm(z_0)}$  for
 $z_0=0.7+0.65 i$.}
\end{figure}

In complete analogy to the standard euclidean Archimedean spirals 
\begin{eqnarray*}
& & \Seucp(w_0)=\left\{|w|
e^{i \varphi} \in \C \, : \,  |w|=+\varphi+|w_0|-\arg w_0, \, \varphi \ge
\arg w_0-|w_0| \right\}\,
, \\[1mm]
& & \Seucm(w_0)=\left\{|w|
e^{i \varphi} \in \C \, : \,  |w|=-\varphi+|w_0|+\arg w_0, \, \varphi \le
\arg w_0+|w_0| \right\}\,
,
\end{eqnarray*}
which pass through both the origin and the point $w_0 \in \C$,
 we call the curves $\Shypp(z_0)$ and $\Shypm(z_0)$ the
\textit{(standard) hyperbolic Archimedean spirals} through the origin
and the point $z_0 \in \D$.

\medskip

Now fix $z_0 \in \Shyppm(z_0)$ and move towards the origin while staying on 
$\Shyppm(z_0)$. Stop either when you reach the point $z_1^{\pm}$ of 
hyperbolic distance $\dD(0,z_0)-\pi$ from the origin or when you reach the origin.
In the first case, we see from the definition of $\Shyppm(z_0)$
that $\arg z_1^{\pm}=\mp \pi +\arg z_0 $, so $z_1:=z_1^{+}=z_1^{-}.$ 
 Let $z_1:=0$ in the second case.  
In both cases we define
$$\gamma^{\pm}(z_0):=\text{ connected part of } \Shyppm(z_0) \text{ between } z_0
\text{ and } z_1\, . $$

\vspace*{-0.4cm}
\begin{figure}[h]
\centerline{\includegraphics[width=4cm]{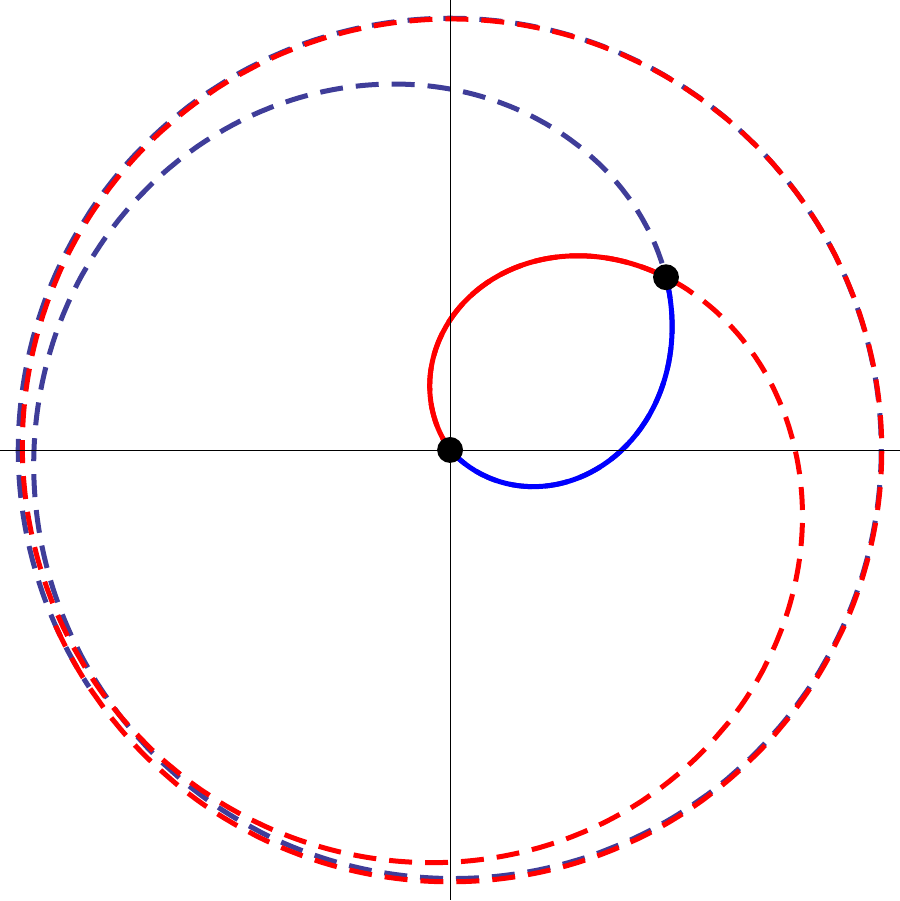}\hspace{2cm}
\includegraphics[width=4cm]{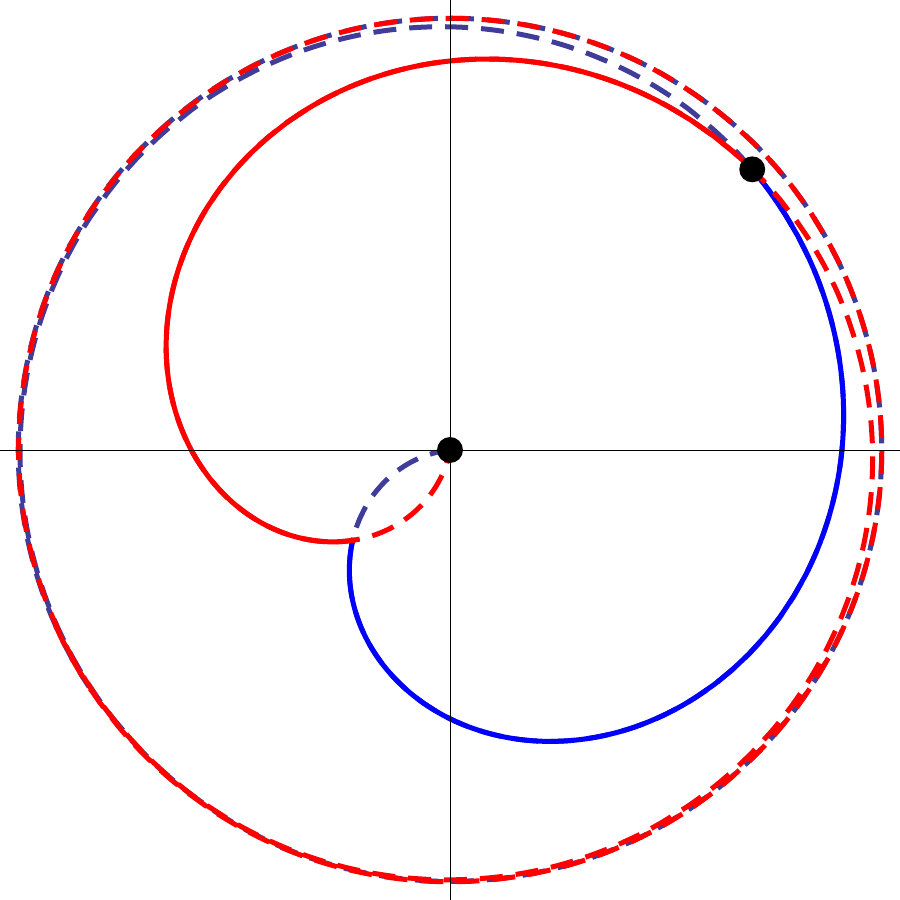}}
\caption{$\Shyppm(z_0)$ (dashed) and $\gamma^{\pm}(z_0)$ (solid) in blue ($+$) and
in red ($-$) for $z_0=0.5+0.4 i$ (left) and for $z_0=0.7+0.65 i$ (right); }
\end{figure}

\begin{corollary} \label{thm:2}
Let $z_0 \in \D$. Then the region of values $\mathcal{V}(z_0)$ has the
following properties.
\begin{itemize}
\item[(a)]
$ \partial \mathcal{V}(z_0)=\gamma^+(z_0) \cup \gamma^-(z_0) \cup \{ 0\}$ and
$\mathcal{V}(z_0) \cup \{ 0\}$ is a Jordan domain bounded by the Jordan curve
$\gamma^+(z_0) \cup \gamma^-(z_0)$.
\item[(b)] The origin is an isolated boundary point of $\mathcal{V}(z_0)$ if and only
  if $|z_0| > \tanh (\pi/2) =0.917152\ldots$\,.
\item[(c)] $\mathcal{V}(z_0)$ is convex if and only if $|z_0| \le
  \tanh(\pi/4)=0.655794\ldots$\,.\footnote{The number $\tanh(\pi/4)$ is also the radius
    of starlikeness in the class $\mathcal{S}$ found by Grunsky 1934.}
\item[(d)] Each boundary curve $\gamma^{\pm}(z_0)$ has hyperbolic length 
$$ L_h(\gamma^{\pm}(z_0))=\begin{cases} \sinh(\dD(0,z_0)) & \text{ if } |z_0| \le
  \tanh(\pi/2)\, , \\
\sinh(\dD(0,z_0))-\sinh(\dD(0,z_0)-\pi)  & \text{ if } |z_0| \ge \tanh(\pi/2) \, .\\
\end{cases}
$$
\end{itemize}
\end{corollary}

Theorem \ref{thm:1} as well as Corollary \ref{thm:2} will be proved
in Section \ref{sec:proof} by making use of the Loewner differential equation
\begin{equation} \label{eq:loew}
\begin{array}{rcl}
\dot{w}(t)&=&-w(t) \, \displaystyle \frac{\kappa(t)+w(t)}{\kappa(t)-w(t)} \, , \qquad t \ge 0 \,
, \\[3mm]
w(0) &=& z_0 \in \D \, ,
\end{array}
\end{equation}
with \textit{fixed} $z_0 \in \D$, 
and its \textit{reachable set} 
$$ \begin{array}{rcl}
\mathcal{R}(z_0) := \{ w(t)  &: &  w(\cdot) \text{ is a solution to
  (\ref{eq:loew}) for some}\\ & & \text{ continuous function }  \kappa : [0,\infty)
\to \partial \D , \, t \ge 0 \}\, , \qquad z_0 \in \D 
\, .
\end{array}$$

\hide{It turns out that the following result holds.

\begin{theorem} \label{thm:2b}
$\mathcal{V}(z_0)=\mathcal{R}(z_0)$ for every $z_0 \in \D$.
\end{theorem}}

Note that by standard results of Loewner's theory it is clear that
 $\overline{\mathcal{V}(z_0)}=\overline{\mathcal{R}(z_0)}$.
It will turn out that actually the stronger statement
$\mathcal{V}(z_0)=\mathcal{R}(z_0)$ holds, so in order to prove Theorem
\ref{thm:1} it suffices to describe the reachable set $\mathcal{R}(z_0)$. By definition,  $\mathcal{R}(z_0)$ is made up of 
the trajectories of the Loewner equation (\ref{eq:loew}).
Any trajectory $w : [0,T] \to \D$ of the Loewner equation 
such that
$w(T)$ is a boundary point of $\mathcal{R}(z_0)$ is therefore of special
interest and
is called an \textit{optimal
trajectory}. The corresponding ``control function'' $\kappa : [0,T]
\to \partial \D$ is then called optimal.

\begin{remark}[The principle of optimality]
The well--known principle of optimality asserts that if a control function
$\kappa : [0,T] \to \D$ ist optimal, then it is also optimal on $[0,T']$
for any $T' \le T$. This means that a trajectory $w : [0,T] \to \D$ of
the Loewner equation which delivers at time $t=T$ a boundary point $w(T)$ of
the reachable set $\mathcal{R}(z_0)$, then 
 $w(t) \in \partial \mathcal{R}(z_0)$ for \textit{every} $t \in [0,T]$.
Hence every initial piece of an optimal trajectory is again optimal.
This is a general principle which holds for any control system under very
general assumptions. It holds in particular for the Loewner equation and has
been utilized in this context before e.g.~by Friedland and Schiffer \cite{FS1977},
Kirwan and Pell \cite{KS1982} and more recently by Graham, Hamada, Kohr and
Kohr \cite{GHKK} for the Loewner equation in several complex variables.
\end{remark}

The next result shows that for any $z_0 \in \D$
the Loewner equation (\ref{eq:loew}) has exactly two optimal trajectories. 
These optimal trajectories  parametrize the entire boundary
of the reachable set $\mathcal{R}(z_0)$. They form exactly the arcs $\gamma^{\pm}(z_0)$
 of the  hyperbolic Archimedean spirals $\Shyppm(z_0)$.

\begin{theorem} \label{thm:3}
Let $z_0 \in \D$. 
Then there exist two trajectories
$w^{\pm}_{z_0} : [0,\infty) \to \D$ such that
$w^{\pm}_{z_0}([0,\infty))$ is the connected part of $\Shyppm(z_0)$ between $z_0$
and $0$.
\begin{itemize}
\item[(a)] If $\dD(0,z_0) \le  \pi$, then $w^{\pm}_{z_0} :[0,\infty) \to \D$ is
  optimal for every $t>0$ and $w^{\pm}_{z_0} ([0,\infty))=\gamma^{\pm}(z_0)
  \backslash \{0\}$.
\item[(b)] If $\dD(0,z_0) >  \pi$, then $w^{\pm}_{z_0} :[0,T] \to \D$ is
  optimal for every $0<T\le T_{max}(z_0)$, but not optimal for
  any $T>T_{max}(z_0)$. Here
$$ T_{max}(z_0)=-\log \left[ 
    \frac{\sinh(\dD(0,z_0)-\pi)}{\sinh(\dD(0,z_0))} \right] \, .$$
Furthermore, $w^{\pm}_{z_0}([0,T_{max}(z_0)])=\gamma^{\pm}(z_0)$.
\end{itemize}

\end{theorem}

\begin{remark}
Theorem \ref{thm:1} has a very well--known classical counterpart for the class
$$\mathcal{S}:=\left\{f : \D \to \C \text{ univalent}, f(0)=0, f'(0)=1\right\}$$
established by Grunsky 1932, see \cite{Grunsky}, who proved the remarkable fact that the set $$\mathcal{W}(z_0):=\left\{\log
(f(z_0)/z_0) \, : \, f \in \mathcal{S}\right\}$$ is exactly the disk
$$ \left\{ w \in \C \, : \, \left|w-\log \frac{1}{1-|z_0|^2} \right| \le
  \log \frac{1+|z_0|}{1-|z_0|} \right\} \, .$$
Grunsky's result was extended e.g.~by Gorjainov and Gutljanski \cite{GG}, who obtained a precise description of
the sets
$$ \mathcal{W}^M(z_0):=\left\{\log \frac{f(z_0)}{z_0} \, : \, f \in \mathcal{S} \text{ such that } |f(z)| < M \text{ for every } z
\in \D\right\}$$
for any $M>1$. The results in \cite{Grunsky, GG} are much more difficult to prove than the results
of the present paper since the sets $\mathcal{W}(z_0)$ and $\mathcal{W}^M(z_0)$
have in fact a much more complicated structure than the set
$\mathcal{V}(z_0)$. In principle, the set $\mathcal{V}(z_0)$ can certainly be described
using the results in \cite{GG} about the sets $\mathcal{W}^M(z_0)$. The purpose of the
present paper is to give a \textit{simple} and \textit{direct} proof of the \textit{simple}
structure of the set $\mathcal{V}(z_0)$ emphasizing some of its remarkable hyperbolic
geometric properties without making appeal to the deeper results about the
sets $\mathcal{W}(z_0)$ and $\mathcal{W}^M(z_0)$ due to Grunsky \cite{Grunsky}
and Gorjainov and Gutljanski \cite{GG}.
\end{remark}

\section{Results: The upper half-plane}

We now replace the interior normalization $f(0)=0$ and $f'(0) \ge 0$, which was
used throughout section \ref{sec:1}, by an appropriate boundary condition. For
this purpose, it is convenient to use the upper half-plane $\H:=\{z \in \C \, :
\, \Im z>0\}$ with its distinguished boundary point $\infty$ instead of the
unit disk $\D$. We  consider the set $ \mathcal{H}_{\infty}(\H)$ of all
holomorphic functions $g : \H \to \H$ such that the so--called hydrodynamic normalization,
$$ g(z)-z \to 0 \text{ as } z\to \infty  \text{ in } S_{\beta}:=\{z \in \C \, : \, |\arg z-\pi/2|<\beta\}  \, , \,
0<\beta<\pi/2 \, 
, $$
is satisfied.

\begin{remark}
The hydrodynamic normalization for functions in  $\mathcal{H}_{\infty}(\H)$
corresponds to the interior normalization for functions in $\mathcal{H}_0(\D)$
in the following way. First note that the interior condition $f(0)=0$ and $f'(0) \ge 0$ for
every $f \in \mathcal{H}_0(\D)$
enforces that the only conformal automorphism of $\D$ contained in $\mathcal{H}_0(\D)$
is the identity. Now, every function $g \in \mathcal{H}_{\infty}(\H)$  has an angular
derivative at $z=\infty$,
$$ g'(\infty):=\angle \lim \limits_{z \to \infty} \frac{g(z)}{z}=1 \, ,$$
and hence the angular limit
$$ g(\infty):=\angle \lim \limits_{z \to \infty} g(z)=\infty\, . $$
Note that the two boundary conditions $g(\infty)=\infty$ and $g'(\infty)=1$
alone allow for infinitely many conformal automorphisms $z+b$, $b \in \R$, of
$\H$, but the sharper hydrodynamic condition for $\mathcal{H}_{\infty}(\H)$ guarantees that the only conformal automorphism of
$\H$ contained in  $\mathcal{H}_0(\H)$
is the identity.
\end{remark}

The following result is the analogue of Theorem \ref{thm:1} for the upper
half-plane.

\begin{theorem} \label{thm:halfplane1}
Let $z_0 \in \H$. Then
$$ \left\{g(z_0) \, : \, g \in \mathcal{H}_{\infty}(\H) \text{ univalent}\right\}=\left\{z
\in \C \, : \, \Im z>\Im z_0\right\} \cup \{z_0\}\, .$$
\end{theorem}

\begin{remark}  \label{rem:halfplane1}
Using the well--known Nevanlinna representation 
for holomorphic functions in $\H$ with positive imaginary part (see
\cite[Theorem 5.3]{RR94}), it is immediate that
$$ \{g(z_0) \, : \, g \in \mathcal{H}_{\infty}(\H)\} \subseteq \{z
\in \C \, : \, \Im z>\Im z_0\} \cup \{z_0\}\, .$$
Hence, Theorem \ref{thm:halfplane1} tells us that the set  of values
$g(z_0)$ for all \textit{univalent} functions $g \in \mathcal{H}_{\infty}(\H)$
is the same as the set of values $g(z_0)$ for all $g \in
\mathcal{H}_{\infty}(\H)$. This is a significant difference to the unit disk
case, where a comparison of Theorem \ref{thm:1} with Rogosinski's Lemma shows that
the set $\mathcal{V}(z_0)$ of values $f(z_0)$ for all univalent
functions $f \in \mathcal{H}_{0}(\D)$ is strictly smaller than the set of values $f(z_0)$ for all functions $f \in \mathcal{H}_{0}(\D)$. 
\end{remark}

As in the unit disk case, the proof of Theorem \ref{thm:halfplane1} relies on
the Loewner differential equation, but now we have to use the \textit{chordal}
version
\begin{equation} \label{eq:chor}
\begin{array}{rcl}
\dot{w}(t)&=& \, \displaystyle \frac{-2}{w(t)-U(t)} \, , \qquad t \ge 0 \,
, \\[3mm]
w(0) &=& z_0 \in \Ha \, ,
\end{array}
\end{equation}
with fixed $z_0 \in \H$, and its reachable set
$$ \begin{array}{rcl}
\mathcal{R}(z_0) := \{ w(t)  &: &  w(\cdot) \text{ is a solution to
  (\ref{eq:chor}) for some}\\ & & \text{ continuous function }  U : [0,\infty)
\to \R , \, t \geq 0 \}\, , \qquad z_0 \in \Ha
\, .
\end{array}$$
It is well--known (see \cite[Chapter 4]{Law2005} and \cite{CDG,GM}) that every solution
$g_t$ of the chordal Loewner equation
$$ \dot{g}_t(z)=\frac{-2}{g_t(z)-U(t)} \, , \qquad g_0(z)=z \, ,$$
generated by a continuous driving function $U : [0,\infty) \to \R$ is a univalent
function on $\H$ such that
$$ g_t(z)=z-\frac{b(t)}{z}+ \mathcal{O}(1/z^2) \quad \text{ as  } z \to \infty \,
.$$
In particular,  $g_t \in \mathcal{H}_{\infty}(\H)$ for every $t \ge 0$. Hence
Remark \ref{rem:halfplane1} implies that
Theorem \ref{thm:halfplane1} is an immediate consequence of the following
result, which describes the reachable set $\mathcal{R}(z_0)$ for any $z_0 \in \H$. 

\begin{theorem} \label{thm:halfplane2} 
Let $z_0 \in \H$. Then
$\mathcal{R}(z_0)=\{z\in\Ha  \, : \, \Im z > \Im z_0 \}\cup
  \{z_0\}$.
\end{theorem}

\section{Proofs: the unit disk} \label{sec:proof}

The proofs are based on an  extremely simple differential inequality 
for the ``hyperbolic polar coordinates'' of the solutions of Loewner's equation,
which follows immediately from the particular form of the Loewner equation.

\subsection{The basic differential inequality}

Fix $z_0 \in \D$ and let
 $w : [0,\infty) \to \D$ be the solution to the Loewner equation
(\ref{eq:loew}) generated by a measurable function $\kappa : [0,\infty) \to \partial \D$, that is,
$$ \frac{d}{dt} \big\{ \log w(t)\big\}
=-\frac{\kappa(t)+w(t)}{\kappa(t)-w(t)} \, , \qquad w(0)=z_0 \, .$$
Separation into real and imaginary parts and writing $w(t)=|w(t)| e^{i
  \varphi(t)}$ gives the equivalent pair of equations
\begin{align} \label{eq:loew1a}
 \frac{d}{dt} \big\{ |w(t)| \big\} &=-|w(t)|
\frac{1-|w(t)|^2}{|\kappa(t)-w(t)|^2} 
\intertext{and} \label{eq:loew2}
\frac{d}{dt} \big\{\varphi(t)\big\} &=-2 \, \frac{\Im\left\{\overline{\kappa(t)} \,
    w(t)\right\}}{|\kappa(t)-w(t)|^2} \, .
\end{align}
Using (\ref{eq:hypdist}) we can rewrite (\ref{eq:loew1a}) in the form
\begin{equation} \label{eq:loew1}
\frac{d}{dt} \big\{\dD(0,w(t))\big\}= \frac{d}{dt} \left\{ \log \frac{1+|w(t)|}{1-|w(t)|} \right\} =
-2 \, \frac{|w(t)|}{|\kappa(t)-w(t)|^2}. \, 
\end{equation}

Now using the simple inequality
\begin{equation} \label{eq:q1}
 |w(t)| \ge \pm \Im \left\{ \overline{\kappa(t)} \, w(t) \right\} \, ,
\end{equation}
 we combine (\ref{eq:loew2}) and (\ref{eq:loew1}) and 
arrive at the differential inequality
\begin{equation} \label{eq:in}
 \frac{d}{dt} \big\{ \dD(0,w(t) \big\} \le \pm
 \frac{d}{dt} \big\{ \varphi(t) \big\} \, .
\end{equation}
\hide{By (\ref{eq:loew1}), the left hand side of (\ref{eq:in}) never vanishes,
so also the right  side of (\ref{eq:in}) is never zero. In view of (\ref{eq:loew2}) 
this means that  $\Im \left\{\overline{\kappa(t)} \,
  w(t) \right\}$ never vanishes.
 Therefore $t \mapsto
\varphi(t)$ is either strictly monotonically decreasing (Case a)
or strictly monotonically increasing (Case b), so
either
\begin{align} \label{eq:ina} \tag{2.4a}
\frac{d}{dt} \left\{ \log \frac{1+|w(t)|}{1-|w(t)|} \right\} & \le +
 \frac{d}{dt} \big\{ \varphi(t) \big\} \, \\
\intertext{or}  
\label{eq:inb}  \tag{2.4b}
\frac{d}{dt} \left\{ \log \frac{1+|w(t)|}{1-|w(t)|} \right\} & \le -
 \frac{d}{dt} \big\{ \varphi(t) \big\} \, .
\end{align}}

\subsection{Integrating the basic differential inequalities and proof of
  Theorem \ref{thm:1}}

Integrating (\ref{eq:in})  yields
\begin{equation} \label{eq:set}
  \dD(0,w(T))-\dD(0,z_0) \le -\left| \varphi(T)-\arg z_0 \right| \, , \qquad T
  \in [0,\infty) \, ,
\end{equation}

with equality  if
\begin{align} \label{eq:ina} \tag{Case I}
 \frac{d}{dt} \big\{ \dD(0,w(t)) \big\}&=+\frac{d}{dt}
 \big\{\varphi(t)\big\} \le 0 \, \,\qquad \text{ for all }  t \in [0,T]\\ 
\intertext{or if} \label{eq:inb} \tag{Case II}
\frac{d}{dt} \big\{ \dD(0,w(t))  \big\} &=-\frac{d}{dt}
\big\{\varphi(t)\big\} \le 0  \qquad \text{ for all } t \in [0,T] \,.
\end{align}

This shows that
\begin{equation} \label{eq:incl}
\mathcal{R}(z_0) \subseteq 
\left\{w=|w| e^{i \varphi} \in \D \, : \,
\dD(0,w)-\dD(0,z_0)\le -|\varphi-\arg z_0| , \, \varphi \in \R \right\}\,.
\end{equation}

We now show that this inclusion is sharp.
 We first consider Case I. Then
equality holds in (\ref{eq:set}) in this case, i.e.,
$$ \dD(0,w(T))-\dD(0,z_0) =\varphi(T)-\arg z_0 \, ,$$
and  equality holds in (\ref{eq:q1}) with $+$  for every $t \in [0,T]$.
This means that
\begin{equation} \label{eq:kap}
\overline{\kappa(t)} w(t)=i |w(t)| \qquad \text{ for every } t \in [0,T] \, .
\end{equation}

Therefore, (\ref{eq:loew1}) reduces to
$$ \frac{d}{dt}  \left[ \log \frac{1+|w(t)|}{1-|w(t)|} \right] =
-2 \, \frac{|w(t)|}{1+|w(t)|^2} \, $$
or
$$ \frac{d}{dt} \big\{ \dD(0,w(t)) \big\} =-\tanh \left( \dD(0,w(t)) \right)
\, , $$
and an  integration leads to
$$ \dD\left(0,w(t)\right)=\arcsinh \left( e^{-t} \sinh (\dD(0,z_0)) \right) \, .$$

In a similar way,  (\ref{eq:loew2}) reduces to
$$
\frac{d}{dt} \big\{\varphi(t)\big\} =  -\tanh \left( \dD(0,w(t)) \right) \, , $$
and another integration gives
$$ \varphi(t)= \varphi^+(t):=\text{Arcsinh} \left( e^{-t} \sinh
  \left(\dD(0,z_0)\right)  \right)
-\dD(0,z_0) +\arg z_0 \, .
$$
Using (\ref{eq:kap}), we see that
$\kappa(t)=-i e^{i \varphi^+(t)}$ for every $t\in [0,T]$. 
We can now \textit{define} 
$$ \kappa^+(t):=-i e^{i \varphi^+(t)} \, , \qquad t \in [0,\infty) \, .$$
This gives a (real--analytic) control function $\kappa^+ : [0,\infty)
\to \partial \D$.  The corresponding solution $w^+_{z_0} : [0,\infty) \to \D$
of the Loewner equation (\ref{eq:loew})  
satisfies
\begin{equation} \label{eq:wplus}
 \dD(0,w^+_{z_0}(t))-\dD(0,z_0)=\arg w^+_{z_0}(t)-\arg z_0 \text{ for all } t
\ge 0 \, 
\end{equation}
by construction. Case II can be handled in a similar way. In fact, if we set
$$ \kappa^-(t):=-i e^{i \varphi^-(t)}$$
where
$$\varphi^-(t):=-\arcsinh \left( e^{-t} \sinh\left( \dD(0,z_0) \right)  \right)
+\dD(0,z_0) +\arg z_0 \, ,$$
then the corresponding solution $w^-_{z_0} : [0,\infty) \to \D$
of the Loewner equation (\ref{eq:loew})  satisfies
\begin{equation} \label{eq:wminus}
 \dD(0,w^-_{z_0}(t))-\dD(0,z_0)=-\arg w^-_{z_0}(t)+\arg z_0 \text{ for all } t
\ge 0 \, .
\end{equation}
This shows that the inclusion (\ref{eq:incl}) is sharp. In order to conclude
\begin{equation} \label{eq:loew20}
 \mathcal{R}(z_0)=\left\{w=|w| e^{i \varphi} \in \D \, : \,
\dD(0,w)-\dD(0,z_0)\le -|\varphi-\arg z_0| , \, \varphi \in \R \right\}\backslash \{
0\} \,  ,
\end{equation}
we just note that $\mathcal{V}(z_0) \cup \{0\}$ is a compact subset of the
complex plane, which is starlike with respect to the boundary point $0$ (i.e.,
if $w \in \mathcal{V}(z_0)$, then $rw \in \mathcal{V}(z_0)$ for any $0 \le r
  \le 1$) and that
 $\overline{\mathcal{R}(z_0)}=\overline{\mathcal{V}(z_0)}$ by
Loewner's theory. This proves  (\ref{eq:loew20}) and 
 $\mathcal{R}(z_0)=\mathcal{V}(z_0)$. The proof of Theorem
\ref{thm:1} is complete.

\subsection{Proof of Theorem \ref{thm:3}}
We  note that in view of (\ref{eq:wplus}) and (\ref{eq:wminus})
the two trajectories $w^{\pm}_{z_0}([0,\infty))$ form the
 connected parts  of 
the hyperbolic Archimedean spirals $\Shyppm(z_0)$ between $z_0$
and $0$. These two trajectories never meet  if and only if
$\varphi^+(t)-\arg z_0>-\pi$ and $\varphi^-(t)-\arg z_0<\pi$ for every $t>0$,
which happens if and only if $\dD(0,z_0) \le \pi$.
In this case, we thus have $w^{\pm}_{z_0}([0,\infty))=\gamma^{\pm}\backslash \{ 0\}$.
If $\dD(0,z_0) > \pi$, then  the two trajectories $w^{\pm}_{z_0}$  first meet at time $t=T_{max}(z_0)$, where
$\varphi^+(T_{max}(z_0))=-\pi+\arg z_0$, that is,
$$ \arcsinh \left( e^{-T_{max}(z_0)} \sinh (\dD(0,z_0)) \right)=\dD(0,z_0)-\pi
\, ,$$
i.e., 
$$ T_{max}(z_0)=-\log \left[ \frac{\sinh(\dD(0,z_0)-\pi)}{\sinh(\dD(0,z_0))} \right]\, .$$
Clearly, $w^{\pm}_{z_0}$ is not optimal on $[0,T]$ for any $T >T_{max}(z_0)$
and $w^{\pm}_{z_0}([0,T_{max}(z_0)])=\gamma^{\pm}(z_0)$.
This proves Theorem \ref{thm:3}. 

\subsection{Proof of Corollary \ref{thm:2}}
(a) is clear from the above discussion. (b) The origin is an isolated boundary
point of $\mathcal{V}(z_0)$ if and only if $\dD(0,z_0) >\pi$, that is, $|z_0|
>\tanh(\pi/2)$. (c) $\mathcal{V}(z_0)$ is convex if and only if
$\varphi^+(t)-\arg z_0 \ge -\pi/2$ and $\varphi^-(t)-\arg z_0 \le \pi/2$ for
every $t>0$. This occurs if and only if $\dD(0,z_0) \le \pi/2$, that is, 
$|z_0| \le \tanh(\pi/4)$. 

\smallskip

(d) In hyperbolic polar coordinates the curve
$\gamma^+(z_0)$ has the
parametrization $I\ni \varphi \mapsto |\gamma(\varphi)| e^{i \varphi}$, where
$\dD(0,\gamma(\varphi))=\varphi+\dD(0,z_0)-\arg z_0 $ and 
$$ I:=\begin{cases}
  [\arg z_0-\dD(0,z_0),\arg z_0] & \text{ if } \dD(0,z_0) \le
  \pi \ ,\\
[\arg z_0-\pi,  \arg z_0] & \text{ if } \dD(0,z_0) >
  \pi \ .
\end{cases}$$
Now, a computation shows
$$ 2 \frac{
  |\gamma'(\varphi)|}{1-|\gamma(\varphi)|^2}=\frac{1+|\gamma(\varphi)|^2}{1-|\gamma(\varphi)|^2}=\cosh
\left(\dD(0,\gamma(\varphi))\right) 
\, ,  $$
so
\begin{eqnarray*}
 L_h(\gamma^+(z_0))&=&\int \limits_{I} 2 \frac{
  |\gamma'(\varphi)|}{1-|\gamma(\varphi)|^2} \, d\varphi=\int \limits_{I}
\cosh (\varphi+\dD(0,z_0)-\arg z_0)\, d\varphi\\
&=& \begin{cases} \sinh(\dD(0,z_0)) & \text{ if } \dD(0,z_0) \le \pi \, , \\
\sinh(\dD(0,z_0))-\sinh(\dD(0,z_0)-\pi)  & \text{ if } \dD(0,z_0)>\pi \, .\\
\end{cases}
\end{eqnarray*}

\section{Proofs: the upper half-plane }

Fix $z_0 \in \H$ and let $H(z_0):=\{z \in \H \, : \, \Im z >\Im z_0\}\cup
\{z_0\}$. 
In order to prove Theorem \ref{thm:halfplane2}, we first note that 
for every $g
\in \mathcal{H}_{\infty}(\H)$
\begin{equation} \label{eq:in1}
\Im g(z_0) \ge \Im z_0 \quad \text{ with equality if and only if } g(z)=z \, .
\end{equation}
This follows immediately from the Poisson
representation
$$ \Im g(z)=c \Im z+\frac{\Im z}{\pi} \int \limits_{-\infty}^{\infty}
\frac{d\mu(t)}{(t-\Re z)^2+(\Im z)^2} \, ,$$
where $\mu$ is a nonnegative Borel measure on $\R$ such that
$$\int \limits_{-\infty}^{\infty} \frac{d\mu(t)}{1+t^2}<\infty$$ (see \cite{RR94}) and
$$ c=\lim \limits_{y \to \infty} \Im g(iy)/y\, , $$ i.e., $c=1$ for $f \in \mathcal{H}_{\infty}(\H)$.
Clearly, (\ref{eq:in1}) shows that $\mathcal{R}(z_0) \subseteq H(z_0)$.

\medskip

Now let $z\in H(z_0)$. We need to find a driving function $U$ such that the
solution $w(t)$ of (\ref{eq:chor}) passes through $z$. We separate into real
and imaginary parts and write $w(t)=x(t)+iy(t)$ and $z_0=x_0+iy_0$. Now, we claim
that $U$ can be chosen such that $w(t)$ connects $z_0$ and $z$ by a straight line segment, i.e. 
\begin{equation*} x(t)=c\cdot y(t)+x_0 - c \cdot y_0\, , \end{equation*}
where
$$c=\frac{\Re z-\Re z_0}{\Im z - \Im z_0}\, .$$

In order to prove this, we separate
 equation (\ref{eq:chor}) into real and imaginary parts and obtain
\begin{equation*}
\begin{array}{rcl}
\dot{x}(t)&=& \, \displaystyle \frac{2(U(t)-x(t))}{(U(t)-x(t))^2+y(t)^2} \, ,
\quad  \dot{y}(t)=\, \displaystyle \frac{2y(t)}{(U(t)-x(t))^2+y(t)^2}  \, , 
\end{array}
\end{equation*}
with initial conditions $x(0)=x_0$ and $y(0)=y_0$.
We now assume that $x(t)$ and $y(t)$ are related by
$$U(t)-x(t)=c\cdot y(t) \, .$$ 
Then we get the following initial value problem:
\begin{equation*}
\begin{array}{rcl}
\dot{x}(t)= \, \displaystyle \frac{2c}{(1+c^2)y(t)} \, , \quad  \dot{y}(t)=\, \displaystyle \frac{2}{(1+c^2)y(t)}  \, , \quad
x(0)=x_0, \, y(0)=y_0 \, ,
\end{array}
\end{equation*}
which can be solved directly:
\begin{equation*}
\begin{array}{rcl}
y(t)=\, \displaystyle \sqrt{\frac{4}{1+c^2}t+y_0^2} \qquad \text{and} \qquad
x(t)= \, \displaystyle c y(t)+x_0-cy_0 \, , \quad t \ge 0\,  .
\end{array}
\end{equation*}
Hence if we now \textit{define}
$$ U(t):=c y(t)+x(t)= 2c \sqrt{\frac{4}{1+c^2}t+y_0^2}+ x_0-cy_0 \, ,$$
then by construction the solution $w(t)=x(t)+i y(t)$ of (\ref{eq:chor})
satisfies $x(t)=cy(t)+x_0-c y_0$
 In particular, the trajectory $t \mapsto w(t)$ is the halfline starting at
 $z_0$ through the point $z$, so $z \in \mathcal{R}(z_0)$.
This completes the proof of Theorem \ref{thm:halfplane2}.

\vfill
\hspace{0.9cm}\begin{minipage}{12cm}
Oliver Roth\\
Sebastian Schlei{\ss}inger\\
Department of Mathematics\\
University of W\"urzburg\\
97074 W\"urzburg\\
Germany\\
roth@mathematik.uni-wuerzburg.de\\
sebastian.schleissinger@mathematik.uni-wuerzburg.de
\end{minipage}

\end{document}